\documentclass[11pt]{article}
\usepackage{graphicx}
\input epsf
\usepackage{amssymb}
\newcounter{fig}

\newtheorem{theorem}{Theorem}

\title{Mathematical Table Turning Revisited}
\author{by Bill Baritompa, Rainer L\"owen, Burkard Polster, and Marty Ross}
\begin{document}

\maketitle
\centerline{{\Large Abstract}}
We investigate under which conditions a rectangular table can be placed with all four feet on a ground described by a function $\mathbb R^2\to \mathbb R$. 

We start by considering highly idealized tables that are just simple rectangles. We prove that given any rectangle, any continuous ground and any point on the ground, the rectangle can be positioned such that all its vertices are on the ground and its center is on the vertical through the distinguished point. This is a mathematical  existence result and does not provide a practical way of actually finding a balancing position. 

An old, simple, beautiful, intuitive and applicable, but not very well known argument guarantees that a square table can be balanced on any ground that is not ``too wild'', by turning it on the spot. In the main part of this paper we turn this intuitive argument into a mathematical theorem. More precisely, our theorem deals with rectangular tables  each consisting of a solid rectangle as top and four line segments of equal length as legs.  We prove that if the ground does not rise by more than $\arctan\left (\frac{1}{\sqrt 2}\right) \approx 35.26^\circ$  between any two of its points, and  if the legs of the table are at least half as long as its diagonals, then the table can be balanced anywhere on the ground,
without any part of it digging into the ground, by turning the table on the spot. This significantly improves on related results recently reported on in \cite{Martin} and \cite{Polster1} by also dealing with  tables  that are not square, optimizing the allowable  ``wobblyness'' of the ground, giving  minimal leg lengths that ensure that the table won't run into the ground, and providing (hopefully) a more accessible proof.

Finally, we give a summary of related earlier results, prove a number of related results for tables of shapes other than rectangles, and give some advice on using our results in real life.

\newpage

\section{\%\&\$\@\#!!!}
You sit down at a table and notice that it is wobbling, because it is standing on a surface that is not quite even. What to do? Curse, yes, of course. Apart from that, it seems that the only quick fix to this problem is to wedge something under one of the feet of the table to stabilise it. However, there is another simple approach to solving this annoying problem. Just turn the table on the spot! More often than not, you will find a position in which all four legs of the table are touching the ground. This may seem somewhat counterintuitive. So, why and under what conditions does this trick work?

\section{Balancing Mathematical Tables---a Matter of Existence}

In the mathematical analysis of the problem, we will first assume that the ground is the graph of a function $g:\mathbb R^2\to \mathbb R$, and that a {\em mathematical table} consists of the four vertices of a rectangle of diameter 2 whose center is on the $z$-axis. What we are then interested in is determining for which choices of the function $g$ can a mathematical table be {\em balanced locally}: that is, when can a table be moved such that its center remains on the $z$-axis, and all its vertices end up on the ground. 

We first observe that it is not always possible to balance a mathematical table locally. Consider, for example, the reflectively symmetric function of the angle $\theta$ about the $z$-axis with 
$$
g(\theta)=
\left \{ 
        \begin{array}{l} 2 \quad \mbox{ if } 0\leq \theta < \frac{\pi}{2} \mbox{ or } \pi \leq \theta < \frac{3\pi}{2},  \\ 
        1 \quad  \mbox{ otherwise}. \end{array} \right.
$$
\begin{figure}[h]
\centerline{\epsfbox{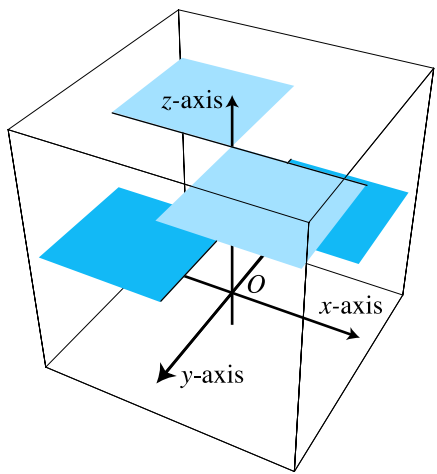}}
\caption[f1]{
 On this discontinuous ground a square mathematical table cannot be balanced locally.\label{cliff}}
\end{figure}

So, the ground consists of four quadrants, two at height 1 and two at height 2; see Figure \ref{cliff}. It is not hard to see that a square mathematical table cannot be balanced locally on such a clifflike piece of ground. On the other hand, we can prove the following theorem:
\begin{theorem}[Balancing Mathematical Tables] A mathematical table can always be balanced locally, as long as the ground function $g$ is continuous.
\end{theorem}

This result is a seemingly undocumented corollary of a theorem by Livesay \cite{Livesay}, which can be phrased as follows: {\em For any continuous function $f$ defined on the unit sphere, we can position a given mathematical table with all its vertices on the sphere such that $f$ takes on the same value at all four vertices.} Note that since our mathematical table has diagonals of length 2, its four vertices will be on the unit sphere iff the centers of the table and the sphere coincide. Choose as the continuous function the {\em vertical distance} from the ground,
$$f:\mathbb S^2 \to \mathbb R:(x,y,z) \mapsto z-g(x,y).$$
Note that here and in everything that follows the vertical distance of a point in space from the ground is really a signed vertical distance; depending on whether the point is above, on or below the ground its vertical distance is positive, zero or negative, respectively. Now, we are guaranteed a position of our rectangle with center at the origin such that all its vertices are the same vertical distance from the ground. This means that we can balance our mathematical table locally by translating it this distance in the vertical direction. Easy!

\section{Balancing Real Tables...by Turning the Tables} 
\begin{figure}[h]
\centerline{\epsfbox{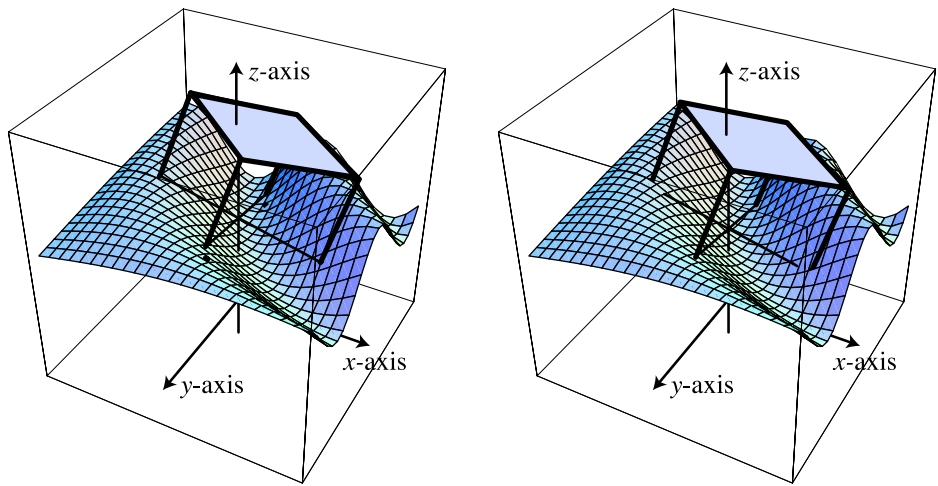}}
\caption[f1]{
 The initial position of the table on the left and its end position on the right. In both positions, two opposite legs are on the ground and the other two are at equal vertical distance from the ground (above the ground in the initial position and below the ground in the end position). Also, in both positions the center of the associated mathematical table is on the $z$-axis. Also check out the Quicktime movie at www.maths.monash.edu.au/$\sim$bpolster/table.mov \label{wobble} to see what happens when we rotate around the $z$-axis.}
\end{figure}

\noindent So, one of our highly idealized tables can be balanced locally on any continuous ground. However, being an existence result, Theorem 1 is less applicable to our real-life balancing act than it appears at first glance. Here are two problems that seem worth pondering:
\begin{enumerate} 
\item {\em Mathematical vs. Real Tables.} A real table consists of four legs and a table top; our theorem only tells us that we can balance the four endpoints of the legs of this real table. However, balancing the whole real table in this position may be physically impossible, since the table top or other parts of the  legs may run into the ground.

To deal with this complication, we  define a  {\em real table} to consist of a solid rectangle with diameters of length 2 as {\em top}, and four line segments of equal length as {\em legs}.  These legs are attached to the top at right angles, as shown in Figure \ref{wobble}. The end points of the legs of a real table form its {\em associated mathematical table}. We say that a real table is balanced locally if its associated mathematical table is balanced locally, and if no point of the real table is below the ground.
\item {\em Balancing by Turning.} A second problem with our analysis so far is that Theorem 1, while guaranteeing a balancing position,  provides no practical method for finding it. After all, although we restrict the center of the table to the $z$-axis, there are still four degrees of freedom to play with when we are actually trying to find a balancing position. 

The following rough argument indicates how, by turning a table on the spot in a certain way, we should be able to locate a balancing position, as long as we are dealing with a square table and a ground that is not ``too crazy''.  
\end{enumerate}

\fbox{\fbox{\parbox{16 cm}{ {\bf Balancing a Square Table by Turning--the Intermediate Value Theorem in Action} \newline Consider a wobbling square table. We wobble the table until two opposite vertices of the associated mathematical table are on the ground, and the other two vertices are the same vertical distance above the ground; see the left diagram in Figure \ref{wobble}. Let's call this position of the table its {\em initial position}. By pushing down on the table, we can make the hovering vertices touch the ground and, in doing so, we have shoved the ``touching'' vertices that same vertical distance into the ground. We call this new position of the table its {\em end position}; see the right diagram in Figure \ref{wobble}. Starting in the initial position, we now rotate the table around the $z$-axis; in doing so, we ensure that at all times the center of the mathematical table is on the $z$-axis, that  the same pair of vertices as in the initial position are touching the ground, and that the other two vertices are an equal vertical distance from the ground. Eventually, we will arrive at the end position. So, we started out with two vertices hovering above the ground, and we finished with the same vertices shoved below the ground. Furthermore, the vertical distance of the hovering vertices depends continuously on the rotation angle. Hence, by the Intermediate Value Theorem, somewhere during the rotation these vertices are also touching the ground: that is, the table has been balanced locally.
}}}

Unlike most other  real-world applications of the Intermediate Value Theorem, it seems that this neat argument is not as well-known as it deserves.  We have not been able to pinpoint its origin, but from personal experience we know that the argument has been around for at least thirty five years and that people keep rediscovering it. In terms of proper references in which variations of the argument explicitly appear, we are only aware of \cite{gardner}, \cite{gardner1}, \cite{gardner2} (Chapter 6, Problem 6), \cite{Hunzinker},  \cite{Kraft}, \cite{Martin}, \cite{Polster}, \cite{Polster1} and \cite{Polster2}; the earliest reference in this list, \cite{gardner}, is Martin Gardner's {\em Mathematical Games} column in the May 1973 issue of {\em Scientific American}. Note that an essential ingredient of the argument is the simple fact that a quarter-turn around its centre takes a square into itself---to move the table from the initial position to the end position takes roughly a quarter-turn around the $z$-axis. Closely related well-documented quarter-turn arguments date back almost a century; see, for example, Emch's proof that any oval contains the vertices of a square in~\cite{Emch} or~\cite{Mayerson}, Section~4. At any rate, we definitely do not claim to have invented this argument.

At first glance, the above argument appears reasonable and, if true, would provide a foolproof method for balancing a square table locally by turning. However, for arbitrary continuous ground functions, it appears just about impossible to turn this intuitive argument into a rigorous proof.  In particular, it seems very difficult  to suitably model the rotating action, so that the vertical distance of the hovering vertices depends continuously upon the rotation angle, and such that we can always be sure to finish in the end position. 

As a second problem, it is easy to construct continuous grounds on which real tables cannot be balanced locally. For example, consider a real square table with short legs, together with a wedge-shaped ground made up of two steep half-planes meeting in a ridge along the $x$-axis. Then it is clear that the solid table top hitting the ground will prevent the table from being balanced locally on this ridge.  

By restricting ourselves to grounds that are not too wild, we can prove that {\em balancing locally by turning} really works. 
\begin{theorem}[Balancing Real Tables] Suppose the ground is described by a Lipschitz continuous function\footnote{Recall that for the ground function $g$ to be Lipschitz continuous means that there
exists a $k$ such that the slope of the line segment connecting any two points on the ground is at 
most $k$: $|g(P)-g(Q)|\leq k|P-Q|$ for any~$P,Q\in {\mathbb R}^2$. The Lipschitz constant of $g$  is then defined to be the  optimal (smallest) choice of $k$. Also, recall that any Lipschitz continuous function is automatically continuous.} 
with Lipschitz constant less than or equal to $\frac{1}{\sqrt 2}$.  Then a real table with ratio $$r=\frac{\mbox {length short side}}{\mbox{length long side}}$$ can be balanced locally on this ground by turning if its legs have length greater than or equal to~$\frac{1}{\sqrt {1+r^2}}$. \end{theorem}

Since $0<r\leq 1$, the maximum of  $\frac{1}{\sqrt {1+r^2}}$ in this range is 1, while all our tables have diagonals of length 2. Thus we conclude that any real table whose legs are at least half as long as its diagonal can be balanced locally by turning on any ``good'' ground. If we are dealing with a square table, then this table can definitely be balanced locally by turning if its legs are at least half as long as its sides. 

Because of the half-turn symmetry of  rectangles, we can be sure to reach a balancing position of a rectangular table whilst turning  it 180 degrees on the spot. As we indicated earlier, to balance a square table we never have to turn it much more than 90 degrees.

For an outline of the following proof for the special case of square tables, aimed at a very general audience, see \cite{Polster1}.  
Furthermore, it has just come to our attention that Andr\'e Martin has also recently published a proof of this result in the special case of square tables and Lipschitz continuous ground functions with Lipschitz constant less than $2-{\sqrt 3}$. In terms of angles, Martin's Lipschitz constant corresponds to 15 degrees and ours, which is  optimal for local turning, to approximately~35.26 degrees.

{\em Proof.} We again start by considering a mathematical table $ABCD$ with diameters of length 2 and centre $O$. Our approach is to bound the wobblyness of our ground by a suitable Lipschitz condition such that putting the two opposite vertices $A$ and $C$ on the ground, and wobbling the table about~$AC$ until $B$ and $D$ are at equal vertical distance from the ground, are unique operations. This ensures that everything in sight moves continuously, as we turn the table on the spot. Following this, it is easy to conclude that we can balance the table locally by turning it. 

Our intuition tells us that to successfully place the four corners, we need four degrees of freedom, four separate motions of the table. Putting our intuition into effect, we approach our balancing act as a succession of four Intermediate Value Theorem (IVT) arguments, taking one ÒdimensionÓ at a time. 

FIRST VERTEX: $A$

Start out with the table hovering horizontally above the ground so that $OA$ lies above the positive $x$-axis, and lower the table until $A$ touches the ground.

SECOND VERTEX: $C$

We now show that since the ground function $g$ is continuous, $A$ can be slid along the ground, towards the $z$-axis (with $O$ sliding up or down the $z$-axis), so that $C$ also touches the ground.
To do this consider the function
$$D(t)= |(t,0,g(t,0))-(-t,0,g(-t,0))|^2=4t^2+(g(t,0)-g(-t,0))^2\, .
$$
Since our table has diagonals of length 2, we want a value of $t\leq \frac{2}{2}=1$ such that $D(t)=2^2=4$. Since $D$ is continuous, $D(0)=0$ and $D(1)\geq 4$, this follows trivially from IVT.

UNIQUENESS OF $C$

Assuming that $g$ is Lipschitz with $Lip(g)\leq 1$, we show that the above positioning of $C$ on the ground
 is unique. This follows from the fact that the  function $D$ is strictly monotonic; this can be seen by differentiating $D$, noting that $|g(t,0)-g(-t,0)|\leq 2t$ with equality only if the function~$g(t,0)$ is linear with slope $\pm 1$ in the interval under consideration. (Lipschitzness is enough for this differentiation argument to work, but a direct algebraic argument is also easy). 

EQUAL HOVERING POSITION

We now rotate the table through an angle $\theta\in [-\frac{\pi}{2}, \frac{\pi}{2}]$ about the diagonal $AC$. We choose the direction so that rotating the table through the angle $-\frac{\pi}{2}$ brings the table into a vertical position with $B$ lying above $AC$. We want to prove the existence of a $\theta$ for which the points $B$ and $D$ are at an equal vertical distance from the ground: we call such a position an {\em equal hovering position}. To show that there is such a special position,  we first choose $\theta=-\frac{\pi}{2}$. The table is now vertical, with~$B$ above $AC$ and $D$ below $AC$.  Since the segments $AB$ and $BC$ are orthogonal, one of the slopes\footnote{We emphasize that  the slope of a line in space is always nonnegative.} of these segments will be greater than or equal to 1. Hence, since $Lip(g)\leq 1$ and since both~$A$ and~$C$ are on the ground, we conclude that $B$ is above or on the ground; similarly, we conclude that $D$ is below or on the ground. If we now rotate the table about $AC$ until $\theta=\frac{\pi}{2}$, then~$B$ is below or on the ground and $D$ is above or on the ground. Now, a straightforward application of IVT guarantees 
a value of $\theta$  for which $B$ and $D$ are an equal vertical distance from the ground. 

UNIQUENESS OF THE EQUAL HOVERING POSITION

We now fix $k\leq 1$ and take the ground to have Lipschitz constant at most $k$. 
We show there exists a choice of
  $k$ which guarantees the  uniqueness of the hovering position. 

Take $A$ and $C$ to be touching the ground as above, with $AC$ then inclined at an angle  $\phi$.
In the following, we sometimes need to express the various objects  as functions of $\theta$, the rotation angle about $AC$ (when assuming the inclination angle $\phi$ to be fixed, which is the case when we are referring to a particular ground); then, for example, $AB$ would be expressed as $AB(\theta)$.  At other
times, we need to express the  objects as functions of  $\phi$ and $\theta$ (when we are not referring to a particular ground); $AB$ would then be expressed as $AB(\phi, \theta)$. 
Here $\phi\in [-\frac{\pi}{4}, \frac{\pi}{4}]$ and $\theta\in [-\frac{\pi}{2}, \frac{\pi}{2}]$. 

We first note that for any equal hovering position the slopes of both $AB$ and
 $BC$ must be at most~$k$ in magnitude. To see this, suppose~$AB$ has slope greater than $k$.  Then, clearly,~$B$ is either above or below the ground. Since $CD$ is parallel to $AB$, it has the same slope as~$AB$; further, if $B$ is higher than $A$, then $D$ is lower than $C$, and vice versa. Therefore, if~$B$ is above the ground, then~$D$ is below the ground, and vice versa. It follows that equal hovering is impossible.

Second, let $tangentB(\theta)$ and $tangentD(\theta)$ be the tangent vectors to the semi-circles swept out by the points $B(\theta)$ and $D(\theta)$, and let $vertB(\theta)$ and
$vertD(\theta)$ be respectively the vertical distances of~$B$ and $D$ to the ground. Note that we have an equal hovering position iff $vertB(\theta)-vertD(\theta)=0$.
It is easy to see that in the $\theta$-interval where  the slope of $tangentB(\theta)$ is greater than or equal to~$k$,  
 then $vertB(\theta)$ is strictly decreasing.   
And also, since $tangentB(\theta)=-tangentD(\theta)$, $vertD(\theta)$  is  strictly increasing in this interval. Thus, $vertB(\theta)-vertD(\theta)$ is strictly decreasing.

Now, let's choose $$k = \min_{\phi,\theta}\,\max\{slopeAB(\phi, \theta), slopeBC(\phi,\theta), 
slopetangentB(\phi,\theta)\}.$$
Of course, $slopeAB(\phi, \theta)$, $slopeBC(\phi,\theta)$, and $slopetangentB(\phi,\theta)$ denote the slopes of $AB$, $BC$ and the  tangent vector at $B$, respectively, and the minimum is taken over all choices of $\phi\in [-\frac{\pi}{4}, \frac{\pi}{4}]$ and $\theta\in [-\frac{\pi}{2}, \frac{\pi}{2}]$. Also, because of compactness, the minimum above is actually achieved.

Given this choice of $k$, we shall show that in the interval where equal hovering is possible the slope of the tangent is at least $k$. So, in this interval, $vertB(\theta)-vertD(\theta)$ is strictly decreasing and thus the equal hovering position must be unique.

We first show that $k=\frac{1}{\sqrt{2}}$. Note that the vectors $AB(\phi, \theta), BC(\phi,\theta),$ and $tangentB(\phi,\theta)$ are mutually orthogonal. If we then write $(0,0,1)$ in terms of this orthogonal frame and take norms, it immediately follows that 
$$1 =\sin^2\beta_1 + \sin^2\beta_2 +\sin^2\beta_3\, ,
$$
where $\beta_1$, $\beta_2$ and $\beta_3$ are the angles the three vectors make with the $xy$-plane. Therefore, at least one of the  $\sin^2\beta_j$ is at least $\frac{1}{{3}}$, and thus the vertical slope 
($=|\tan\beta_j|$) of the corresponding vector
must be at least $\frac{1}{\sqrt{2}}$. It follows that $k \geq  \frac{1}{\sqrt{2}}$.

\begin{figure}[h]
\centerline{\epsfbox{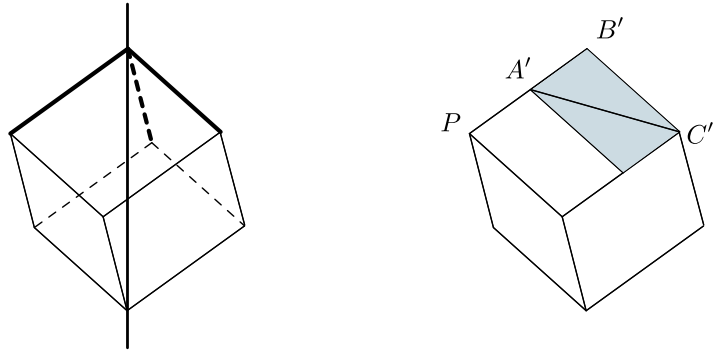}}
\caption[f1]{
 \label{optimal} Three mutually orthogonal segments of length 1 balanced on the $xy$-plane. Then (grey) rectangles of any shape can be fitted as indicated.\label{check}}
\end{figure}

To demonstrate the minimum $k=\frac{1}{\sqrt{2}}$ is achieved, we show any table can be oriented
in the critical position, with all three slopes equal to $\frac{1}{\sqrt{2}}$, 
and with associated tilt angle $\phi$  between~$-\frac{\pi}{4}$ 
and~$\frac{\pi}{4}$.
To do this, consider the tripod formed from three edges
of a cube tilted to have vertical diagonal shown in the left diagram in Figure~\ref{check}. These edges are mutually orthogonal, and one easily 
calculates that the slopes of all three  edges are $\frac{1}{\sqrt{2}}$. Notice that 
every table is similar to one of the grey rectangles, shown in the right diagram, created by moving the point $A'$ from $P$ to~$B'$. Furthermore, it is clear that the slope of the diagonal~$A'C'$ is less than the slope of $B'C'$, which is equal to~$\frac{1}{\sqrt 2}$, guaranteeing the tilt angle $\phi$ is in the desired range.
 By scaling and translating the rectangle suitably, and  relabelling the vertices $A',B',$ and $C'$ as~$A,B$, and $C$, respectively, we arrive at the desired orientation of our table.

It remains to show that $k=\frac{1}{\sqrt 2}$ implies that the slope of the tangent is at least $k$ in the interval where equal hovering is possible. Note that in this interval, $slopeAB$ and $slopeBC$ are at most~$\frac{1}{\sqrt 2}$. Then the equation 
$1 =\sin^2\beta_1 + \sin^2\beta_2 +\sin^2\beta_3\, 
$
 implies that $slopetangentB$ is at least $\frac{1}{\sqrt 2}$.

CONTINUITY OF $A$,  $B$, $C$, AND $D$ 

All of the above calculations were performed with $OA$ projecting to the positive $x$-axis. We now consider rotating the table about the $z$-axis (while of course being willing to tilt the table as we rotate).
So let $\gamma$ be the angle the projection of $OA$ makes with the positive $x$-axis.
By our Lipschitz hypotheis, for any $\gamma$ there is a unique equal hovering position (with the projection of $OA$ making the angle $\gamma$, and tilting the table around $AC$ an angle between $-\frac{\pi}{2}$ and $\frac{\pi}{2}$ in a fixed direction).

We need to show that the positions of the four vertices $A$, $B$, $C$ and $D$, of the equally hovering table are continuous functions of $\gamma$. To do this, consider
a sequence $\gamma_n \to \gamma$, with corresponding corner positions $A_n$, $B_n$, $C_n$, and $D_n$ and $A,B,C,$ and $D$. We want to  show $A_n\to A$,  $B_n\to B$, $C_n\to C$, and $D_n\to D$. By compactness, we can take a subsequence so that the corners converge to something: $A_n\to A^\ast$,  $B_n\to B^\ast$, etc. But, by continuity of everything in sight, $A^\ast$ and $C^\ast$ are touching the ground $B^\ast$ and $D^\ast$ are an equal vertical distance from the ground, all four points are corners of the kind of table we are considering, with the projection of $OA^\ast$ making an angle $\gamma$ with the positive $x$-axis. By uniqueness, we must have $A^\ast=A$, $B^\ast=B$, $C^\ast =C$, and $D^\ast =D$, as desired.

BALANCING POSITION

With the uniqueness of the equal hovering position for a given $\gamma$, and with the continuous dependence of this position of the table upon $\gamma$, we conclude that the distance 
that $B$ and $D$ are hovering above the ground is also a continuous function of $\gamma$. If we are dealing with a square table, we can now finish the proof using IVT one more time, as described in the intuitive table turning argument presented at the beginning of this section. 

For a  general rectangle, let the {\em initial position} be an equal hovering position for  which the $z$-coordinate of the center of the table is a minimum, and let the {\em end position}
be an equal hovering positions for which the $z$-coordinate of the center of the table is a maximum. 
Note that the hovering vertices in the initial position must be on or below the ground: if not, we could create a
lower equal hovering position, contradicting minimality, by pushing vertically down on the table until the hovering vertices touch the ground.  Similarly, in the end position the hovering vertices are
on or  above the ground. Now, IVT can be applied to guarantee that among all the equal hovering positions there is at least one balancing position. 

BALANCING REAL TABLES

To balance a real table of side lengths ratio $r$, we determine a balancing position of the associated mathematical table, as described above. We now show that legs of length at least~$\frac{1}{\sqrt{1+r^2}}$ guarantee that, balanced in this position, none of the points of the real table are below the ground. We give the complete argument for a square table, and then describe how things have to be modified to give the result for arbitrary rectangular tables. In the following, we will refer to the four vertices of the table top  as $A',B',C'$ and $D'$, corresponding to the vertices $A,B,C$ and $D$, respectively.

We first convince ourselves that no matter how long are the legs of our table, no part of a leg of the balanced table will be below the ground. Let's consider the orthogonal tripod consisting of $AB$, $AD$ and the leg at $A$. Since the Lipschitz constant of our ground is at most $\frac{1}{\sqrt 2}$, the slopes of $AB$ and $AD$ are less than or equal to this value; thus, arguing as above, we see that the leg must have slope at least ~$\frac{1}{\sqrt 2}$. This implies that  no  leg of our balanced table will dip below the ground.\footnote{For certain grounds with Lipschitz constant~$\frac{1}{\sqrt 2}$ it is possible that a leg of a balanced table may lie along the ground.}

It remains to choose the length of the legs such that no point of the table top of our balanced table will ever be below the ground. First, fix the length of the legs and consider the inverted {\sc solid} circular cone, whose vertex is one of the vertices of our mathematical table, whose symmetry axis is vertical, and whose slope is~$\frac{1}{\sqrt 2}$. Intersecting this cone with the plane
in which the table top lies  gives a conic section which is either an ellipse, a parabola or a hyperbola.\footnote{Parabolas and hyperbolas can occur since the plane that the table top is contained in can have maximum slope greater than~$\frac{1}{\sqrt 2}$.} Note that since we intersect the plane with a solid cone this conic section will be ``filled in''. We can be sure that a point in this plane is not below the ground if it is contained in the conic section.  Therefore, what we want to show is that the union of the four conic sections associated with the four vertices of our mathematical table contains the whole table top. It is clear that the four conic sections are congruent and that any two of them can be brought into coincidence via a translation. Furthermore, given any point of one of these conic sections, this point
and  the respective points in the other three conic sections form a square that is congruent to our table top.  Finally, since the legs have slope of at least $\frac{1}{\sqrt 2}$, the end point of a leg of our table on the plane is contained in the conic section associated with the other end point of this leg. 

\begin{figure}[h]
\centerline{\epsfbox{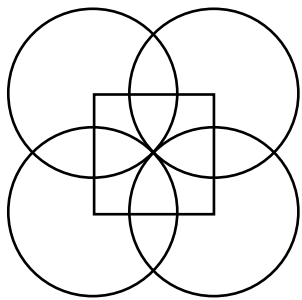}}
\caption[f1]{
 \label{intersect} The four conic sections and the table top in the case of a square table with legs of length~$\frac{1}{\sqrt 2}$ that is balanced horizontally.}
\end{figure}

To show that we need legs of length at least  $\frac{1}{\sqrt 2}$, consider a special ground with Lipschitz constant~$\frac{1}{\sqrt 2}$. This ground coincides with the $xy$-plane outside the unit circle, and above the unit circle it is the surface of the cone with vertex $(0,0,\frac{1}{\sqrt 2})$ and base the unit circle. Since the diagonals of our table are of length 2, the mathematical table will balance locally on this ground iff its vertices are on the unit circle. This means that the length of the legs have to be at least as long as the cone is high if we want to ensure that no point of the table top is below the ground; it follows that we have to choose the length of our legs to be at least $\frac{1}{\sqrt 2}$. If the length of the legs is equal to  $\frac{1}{\sqrt 2}$ and the table is balanced on this ground, then the four conic sections are circles that intersect in the center of the table top as shown in Figure~\ref{intersect}. As you can see, the table top is indeed contained in the union of these four circles. If we make the legs  longer, the circles will overlap more. If we make the legs shorter, the circles will no longer overlap in the middle. 

\begin{figure}[h]
\centerline{\epsfbox{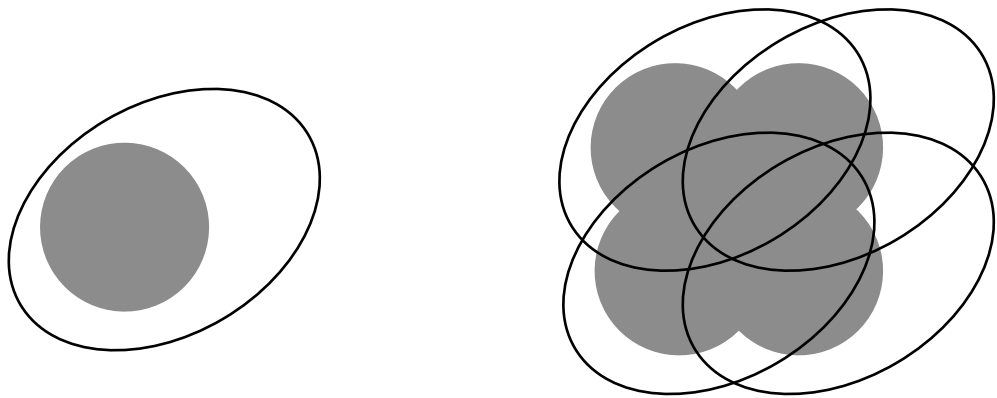}}
\caption[f1]{
 \label{intersect1} If the table is not horizontal, the conic sections are larger than the circular (gray) sections in the horizontal case (left). Their union is simply connected (right).}
\end{figure}

Now, consider any ground, 
and take  the legs to be of length $\frac{1}{\sqrt 2}$; clearly, if we can show that this table does not dip below the ground, then the same is true for any table with longer legs. When we tilt the table away from the horizontal position, the intersection pattern of the conic sections  gets more complicated. 
The critical observation is that tilting the table results in the conic sections getting larger: it can be shown
that each conic
section contains a copy of one of the circles in Figure~\ref{intersect}.\footnote{To see this note that what we are looking at are the possible intersections of a given cone with planes that are a fixed distance from the vertex of the cone.}  
Since, given any point of one of these conic sections, it and  the corresponding points in the other three conic sections form a square that is congruent to our table top, the union of these conic sections will contain a possible translated image of the union of the circles, that we encountered before; see Figure~\ref{intersect1}. So, in a way our previous picture has just grown a little bit and been translated. (Note, however,
it is not immediate that the conic sections together cover the table top rather than some translation of the table
top).
Using this fact and the simple possible convex shapes of the conic sections that we are dealing with, we can conclude that no matter how we tilt the table, the union of these conic sections will always be a simply connected domain. This means that we can be sure that the table top is contained in this union if we can show that the boundary of the table top is contained in it.

\begin{figure}[h]
\centerline{\epsfbox{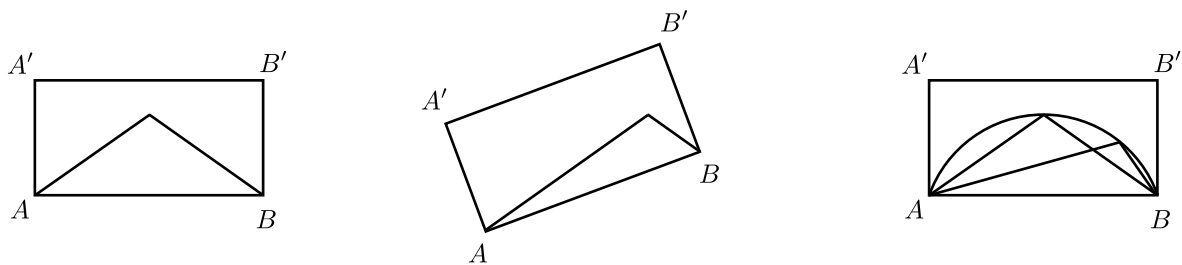}}
\caption[f1]{
 \label{thales} The points of intersection of the two lines of slope $\frac{1}{\sqrt 2}$ for the different possible vertical positionings of the rectangle $AA'BB' $ are on a circle. }
\end{figure}

We proceed to show that for all possible positions of our table in space the sides of the table top never dip below the ground. Clearly, it suffices to show this for one of the sides of the table top, say~$A'B'$.  For this we consider the possible positions of the rectangle with vertices $A, A',B$ and~$B'$ in space. We start with the rectangle vertical and  $AB$ horizontal.  Draw lines of slope~$\frac{1}{\sqrt 2}$ ending in~$A$ and $B$; see Figure~\ref{thales} (left). Since the point of intersection of these two lines is not above~$A'B'$, no point of this segment can be below the ground when the rectangle is positioned in such a  way.  Now rotate the rectangle around its center, keeping it in a vertical plane,  and keeping 
the slope of ~$AB$  less than or equal to~$\frac{1}{\sqrt 2}$. Again, draw lines of slope $\frac{1}{\sqrt 2}$ ending in~$A$ and $B$; see Figure~\ref{thales} (middle). Again, the position of the point at which these two lines intersect tells you whether~$A'B'$ can possibly touch the ground with the rectangle in this position. Since the two lines always intersect in the same angle, we know that the points of intersection are on a circle segment through~$A$ and~$B$; see Figure~\ref{thales} (right).

\begin{figure}[h]
\centerline{\epsfbox{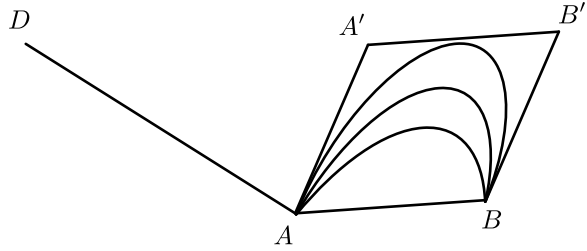}}
\caption[f1]{
 \label{tiltrec} The more we tilt, the closer the apex of the circle segment gets to $A'B'$. As long as the apex is below or on $A'B'$, we can be sure that $A'B'$ does not dip below the ground.}
\end{figure}

Now, tilt the original rectangle around $AB$. We repeat everything that we have done so far to  end up with another circle segment. However, the apex of this circle segment will the  be closer to $A'B'$ than the one we encountered before. In fact, the more we tilt, the closer we will get; see Figure~\ref{tiltrec}.

  Since the apex of one of these circles is the point closest to $A'B'$,  and since  the apex corresponds to~$AB$ being horizontal, we now calculate just how close this apex gets when we tilt around $AB$ though a maximal possible angle. This maximal possible angle is attained if $AD$ (which is orthogonal to~$AB$) has slope $\frac{1}{\sqrt 2}$. It is a routine exercise to check that in this position the slope of the line connecting $A$ with the midpoint of $A'B'$  is $\frac{1}{\sqrt 2}$. This means that  in this position the apex will be contained in $A'B'$. We conclude that if we choose the legs of our square table to be at least~$\frac{1}{\sqrt 2}$ long, then the boundary of the table top, and hence also the table top itself, will not dip below the ground.

For  tables that are not square, the same arguments apply up the point where we start tilting the  rectangle $AA'B'B$ around~$AB$. We now have to worry about two different rectangles corresponding to the longer and shorter sides of the table top. If $r$ is the ratio of the lengths of the sides of the table, then it is easy to see that the critical length of the legs that we need to avoid running into the ground  is the length that makes the longer of the two rectangles similar to the rectangle that we considered in the square case. This critical length is $\frac{1}{\sqrt{ 1+r^2}}$.\hfill $\Box$

\section{Other Balancing Acts}
\subsection*{Horizontal balancing} When we balance a table locally, the table will usually not end up horizontal, and a beer mug placed on the table may still be in danger of sliding off. It would be great if we could arrange it so that the table is not only balanced but also horizontal, maybe by moving the center of the table off the $z$-axis and balancing it somewhere else on the ground. Just imagine the ground to be a tilted plane, and you can see that this will not be possible in general. However, Fenn \cite{Fenn} proved the following result: {\em If a continuous ground coincides with the $xy$-plane outside a compact convex disc and if the ground never dips below the $xy$-plane inside the disk, then a given square table can be balanced  horizontally such that the center of the table lies above the disk.} Let's call the special kind of ground described here a {\em Fenn ground} and the part of this ground inside the distinguished compact disk its {\em hill}. 

The problem of horizontally balancing tables consisting of plane shapes other than squares on Fenn grounds has also been considered. Here `horizontal balancing on a Fenn ground' means that in the balancing position some interior points of the shape are situated above the hill. It has been shown by Zaks  \cite{Zaks} that a triangular table can be balanced on any Fenn ground. In fact, he showed that if we start out with a horizontal triangle somewhere in space and mark a point inside the triangle, then we can balance this triangle on any Fenn ground, with the marked point above the hill, by just translating the triangle. Fenn also showed that tables with four legs that are not concircular and those forming regular polygons with more than four legs cannot always be balanced horizontally on Fenn grounds. Zaks mentions an unpublished proof by L.M. Sonneborn that any polygon table with more than four legs cannot always be balanced horizontally on Fenn grounds. It is not known whether any concircular quadrilateral tables other than squares can always be balanced horizontally  on Fenn grounds. See \cite{Kronheimer}, \cite{Mayerson}, \cite{Mayerson1}, and \cite{Mayerson2} for further results relating to this line of research. 
 
\subsection*{Local Balancing of Exotic Mathematical Tables} Taking things  to different mathematical extreme, we can consider a table consisting of $n\geq 3$ {\em leg points} in 3-space together with an additional {\em center point}. We then ask whether, given any continuous ground, it is possible to always balance this table locally, that is, move this configuration of $n+1$ points into a position in which the $n$ leg points are on the ground, and the center is on the $z$-axis.  

The example of a plane ground shows that the leg points of an always locally balancing table have to be coplanar. Let's consider the example of a ground that contains part of a sphere that is large enough to ensure that all legs of our  table end up on this part of the sphere whenever the table is locally balanced on this ground. Then intersecting this sphere with the plane that the leg points are contained in gives a circle that all leg points are contained in. Hence the leg points of the  table are concircular. Now, let's consider a ground that includes part of an ellipsoid which does not contain a copy of the circumcircle of the leg points of the table; moreover, we choose the ellipsoid large enough so that all leg points of our table end up on the ellipsoid whenever the table is locally balanced on this ground.  Then intersecting the ellipsoid with the plane containing the leg points gives an ellipse that is different from the circumcircle of the leg points. However, this is impossible if the table contains more than four leg points because five points on an ellipse determine this ellipse uniquely.  We conclude that an always locally balancing table must  have three or four leg points and that these points are concircular. Note that requiring concircularity in the case of three points is not superfluous since we need to exclude the case of three collinear points.

Livesay's theorem, which made the proof of Theorem 1 so easy, has a counterpart for triangles, due to Floyd \cite{Floyd}. It is a straightforward exercise to apply this result to prove the following theorem: \begin{theorem}[Balancing Triangular Tables] If the ground function is continuous, a triangular table whose three leg points are contained in a sphere around its center can be balanced locally.
\end{theorem}

Of course, one should be able to prove a lot more when it comes to balancing triangular tables! 

In the case that the center and the (three or four) leg points of an always locally balancing table are coplanar, we can say a little bit more about the location of the center point with respect to the leg points. Begin by balancing the table in the $xy$-plane and drawing the circles around the center that contain leg points; if one of the leg points coincides with the center, then also consider this point as one of the circles. Now it is easy to see that there cannot be more than two such circles. Otherwise a ground that coincides with the $xy$-plane inside the third smallest circle and that lies above the plane outside this circle would clearly thwart all local balancing efforts. Therefore, if we want to check whether our favorite set of three or four concircular points is the set of leg points of a locally balancing table, there are usually very few  positions of the center relative to the leg points which need to be considered. 

Perhaps the most natural choice for the center is the center of the circle that the leg points are contained in. As a corollary to the above theorem for triangles, we conclude that a triangular table with this natural choice of center is always locally balancing. In the case of four concircular points with this natural choice of center we do not know whether any tables apart form the rectangular ones are always locally balancing. However, a result worth mentioning in this context is Theorem~3 in Mayerson's paper~\cite{Mayerson}  (see also the concluding remarks in Martin's paper \cite{Martin}). It can be phrased as follows: {\em Given a continuous ground and one of these special four-legged tables in the $xy$-plane, the table can be rotated in the $xy$-plane around its center such that in this new position the four points on the ground above the leg points are coplanar.} The quadrilateral formed by the coplanar points on the ground will be congruent to the table if and only if the plane it is contained in is horizontal, in which case we have actually found a balancing position for our table. In all other cases, the quadrilateral on the ground is a deformed version of the table. Still, if the ground is not too wild, both quadrilaterals will be very similar, and lifting the table up onto the ground should result in the table not wobbling too much.

Livesay's theorem is a generalization of a theorem by Dyson \cite{Dyson}, which only deals with the square case. A higher-dimensional counterpart of Dyson's theorem arises as a special case of  results of 
Joshi~\cite{Joshi}, Theorem 2 and Yang~\cite{Yang}, Theorem 3: {\em Given a continuous real-valued function defined on the $n$-sphere, there are~$n$ mutually orthogonal diameters of this sphere such that the function takes on the same value at all~$2n$ end points of these diameters.} Note that the endpoints of  $n$ mutually orthogonal diameters of the $n$-sphere are the vertices an $n$-dimensional orthoplex, one of the regular solids in $n$-dimensions. (For example, a 1-dimensional orthoplex is just a line segment and a 3-dimensional orthoplex is an octahedron.) Using the same simple argument as in the case of Livesay's theorem, we can prove the following theorem:
\begin{theorem}[Balancing Orthoplex-Shaped Tables]  An $(n-1)$-dimensional orthoplex-shaped table in $\mathbb R^n$ can be balanced locally on any ground given by a continuous function $g:\mathbb R^{n-1}\to \mathbb R$.
\end{theorem}

For other closely related results see  \cite{de Mira}, \cite{Fenn1}, \cite{Hadwiger}, \cite{Yamabe}, \cite{Yang1}, \cite{Yang2}, and~\cite{Yang3}.

\subsection*{Balance Everywhere} Imagine a square table with diameter of length 2 suspended horizontally high above some ground, with its center on the $z$-axis. Rotate it a certain angle about the $z$-axis, release it, and let it drop to the ground. It is easy to identify continuous grounds such that all four leg of the table will hit the ground simultaneously, no matter what release angle you choose. Of course, any horizontal plane will do, and so will any ground that contains a vertical translate of the unit circle. We leave it as an exercise for the reader to construct a ground that is not of this type but admits horizontal balancing for any angle. Also, the reader may wish to convince themselves that the following is true: we are dealing with a ground as in Theorem 2. If the center of the table has the same $z$-coordinate in all its equal hovering positions (positions in which $A$ and $C$ touch the ground and~$B$ and $D$ are at equal vertical distance from the ground), then in fact the table is balanced in all these positions. 

\section{Some Practical Advice}

\subsection*{Short Legs and Tiled Floors} Note that if you shorten one of the legs of a real-life square table, this table will wobble if you set it down on the plane, and no turning or tilting will fix this problem. In real life rectangular tables the ends of whose legs do not form a perfect rectangle are not uncommon and, as our simple example shows, those uneven legs may conspire to make our anti-wobble tactics fail. 

Considering our example of a discontinuous ground at the beginning of this article, it should be clear that a wobbling table on a tiled floor may also defy our table turning efforts. 

\subsection*{How to Turn Tables in Practice} In practice, it does not seem to matter how exactly you turn your table on the spot, as long as you turn roughly around the center of the table. Notice that you
needn't actually establish the equal hovering: as you rotate towards the correct balancing position,
 there will be less and less wobble-room until, at the correct rotation, the balancing position
 is forced. 
  With a square table, you can even go for a little bit of a journey, sliding the table around in your (continuous) backyard. As long as you aim to get back to your starting position, incorporating a quarter turn in your overall movement, you can expect to find a balancing position.

\end{document}